\documentclass{amsart}
\pdfoutput=1

\usepackage{url}
\usepackage[pdftex,colorlinks,citecolor=black,linkcolor=black,urlcolor=black,bookmarks=false]{hyperref}

\usepackage{doi}
\usepackage{booktabs}

\theoremstyle{plain}
\newtheorem{theorem}{Theorem}

\newtheorem{lemma}[theorem]{Lemma}
\newtheorem{proposition}[theorem]{Proposition}

\numberwithin{theorem}{section}

\numberwithin{equation}{section}

\newcommand{\Z}{{\mathbb Z}}

\title{Chudnovsky's formula for $1/\pi$ revisited}

\begin{document}

\begin{abstract}
The document contains an outline of a modular proof for Ramanujan-Chudnovsky identity
\begin{displaymath}
\frac{1}{\pi}=\frac{12}{640320^{3/2}}\sum_{m=0}^{\infty}(545140134m+13591409)\frac{(-1)^m(6m)!}{640320^{3m}(3m)!(m!)^3}.
\end{displaymath}

\end{abstract}

\author{Yue Zhao}
\address{ECSE Department, Rensselaer Polytechnic Institute, 110 8th Street\\
Troy, NY 12180\\}
\email{yue.zhao756@gmail.com}

\maketitle

\section{Introduction}

Ramanujan listed 14 formulae of $1/\pi$ in his 1914 paper \cite{Ram}, where all the formulae are of the form \begin{equation}\label{genfor}\frac{1}{\pi}=\sum_{n=0}^{\infty}(a+bn)d_nc^n.\end{equation} Two among them are quite impressive:

\begin{equation}\label{rampi1}
\frac{1}{\pi}=\frac{2\sqrt{2}}{99^2}\sum_{m=0}^{\infty}(26390m+1103)\frac{(4m)!}{396^{4m}(m!)^4},
\end{equation}

and

\begin{equation}\label{rampi2}
\frac{1}{\pi}=\frac{2}{84^2}\sum_{m=0}^{\infty}(21460m+1123)\frac{(-1)^m(4m)!}{(84\sqrt{2})^{4m}(m!)^4}.
\end{equation}

What Ramanujan discovered can be seen as examples of level 1 and 2 Ramanujan-Sato series\cite{Chan}. In 1989, David Chudnovsky and Gregory Chudnovsky published\cite{Chudnovsky} a Ramanujan-Sato series of level 1 which converges to $1/\pi$ extremely rapidly:
\begin{equation}\label{chud}
\frac{1}{\pi}=\frac{12}{640320^{3/2}}\sum_{m=0}^{\infty}(545140134m+13591409)\frac{(-1)^m(6m)!}{640320^{3m}(3m)!(m!)^3}.
\end{equation}

We are aiming at giving a (modular) proof of this amazing identity in this document. Our calculation still starts with one of the equalities in \cite{BorBro}, p. 181, formula (5.5.9):

\begin{equation}\label{hyp}({1-4(2kk^\prime)^2})^{1/2}\left(\frac{2K(k)}{\pi}\right)^2={}_3 F_2\left(\begin{matrix}1/6& 5/6& 1/2\\1& 1\end{matrix};-\frac{27(2kk^\prime)^2}{(1-4(2kk^\prime)^2)^3}\right),\end{equation}
 where $k^\prime=\sqrt{1-k^2}$. We can rewrite the formula (\ref{hyp}) as

\begin{equation}\label{Clausen}\left(\frac{2K(k)}{\pi}\right)^2=a(k)\sum_{n=0}^{\infty}b_nc^n(k),\end{equation}
where $K(k)$ is complete elliptic integral of the first kind and $a(k)$, $c(k)$ are rational functions of $k$.

Let $\theta_2, \theta_3, \theta_4$ be the Jacobi theta functions

\begin{displaymath}\theta_2(q)=\sum_{n\in\Z}q^{(n+1/2)^2},\end{displaymath}

\begin{displaymath}\theta_3(q)=\sum_{n\in\Z}q^{n^2},\end{displaymath}

\begin{displaymath}\theta_4(q)=\sum_{n\in\Z}(-1)^nq^{n^2},\end{displaymath}

$\eta$ be the Dedekind $\eta$-function

\begin{displaymath}\eta(q)=q^{1/24}\prod_{n=1}^{\infty}(1-q^n).\end{displaymath} Classical elliptic function theory gives \cite{BorBro}(p. 69)

\begin{equation}\label{Kexp}\left(\frac{2K(k)}{\pi}\right)^2=2^{4/3}\eta^{4}(q^2)(kk^{\prime})^{-2/3}=\theta_3^4(q),\end{equation}

where

\begin{displaymath}k=\frac{\theta_2^2(q)}{\theta_3^2(q)}, k^{\prime}=\frac{\theta_4^2(q)}{\theta_3^2(q)},q=\exp{(-\pi\tau)}.\end{displaymath}

Taking logarithmic differentiation by $k$ on both sides of (\ref{Clausen}) and using formula (2.3.10) in \cite{BorBro},

\begin{displaymath}\frac{\mathrm{d} q}{\mathrm{d} k}=\frac{\pi ^2q}{2k {k^{\prime}}^2K^2},\end{displaymath}

we get

\begin{equation}\label{Clausen-derived}\frac{1}{6}P(q)=u(k)\left(\frac{2K(k)}{\pi}\right)^2+v(k)\sum_{n=0}^{\infty}nb_nc^n(k),\end{equation}

where

\begin{displaymath}P(q)=1-24\sum_{n=1}^{\infty}\frac{nq^{2n}}{1-q^{2n}}.\end{displaymath}

$u, v, c$ are all rational functions of $k$.

Taking logarithmic differentiation on both sides of the transformation formula of $\eta$-function

\begin{equation}\label{etatrans}\eta(q^2(1/\tau))=\tau^{1/2}\eta(q^2(\tau))\end{equation}

would lead to

\begin{displaymath}\tau^2 P(q(\tau))+P(q(1/\tau))=6\tau/\pi.\end{displaymath}

Let $\tau=\sqrt{n}$, $n\in\mathbb{N}$. Then

\begin{equation}\label{etatranslog}nP(e^{-\pi\sqrt{n}})+P(e^{-\pi/\sqrt{n}})=6\sqrt{n}/\pi.\end{equation}

Denote $k_1=k(q^n),k_2=k(q)$. The theory of elliptic functions suggests that $K^2(k_1)/K^2(k_2)$ is an algebraic function of $k_1$ and $k_2$. Taking logarithmic differentiation of $K^2(k_1)/K^2(k_2)$ by $k$ and using formula (\ref{Kexp}), it can be shown that \begin{equation}\label{ramfuncs}nP(q^n)-P(q)=\left(\frac{2K(k)}{\pi}\right)^2G_1\end{equation}

where $G_1$ is an algebraic function of $k_1$ and $k_2$.

Take $q=e^{-\pi/\sqrt{n}}$ in (\ref{ramfuncs}). Direct computation with formula (\ref{Clausen-derived}), (\ref{etatranslog}) and (\ref{ramfuncs}) would establish the following result:

\begin{equation}\label{symbolform}\frac{1}{\pi}=\sum_{m=0}^{\infty}(2\sqrt{n}v(k)m+G_0)b_mc^m(k),\end{equation}
where \begin{math}k=k(e^{-\pi\sqrt{n}})\end{math},

\begin{displaymath}b_m=\frac{(6m)!}{12^{3m}(3m)!(m!)^3},\end{displaymath}

\begin{displaymath}2v(k)=\frac{((k^\prime)^2-k^2)(1+8(2kk^{\prime})^2)}{1-4(2kk^{\prime})^2},\end{displaymath}

\begin{displaymath}c(k)=-\frac{27(2kk^\prime)^2}{(1-4(2kk^\prime)^2)^3},\end{displaymath}
and
\begin{displaymath}G_0=\frac{\sqrt{n}}{6}\left(\frac{((k^\prime)^2-k^2)(1+8(2kk^{\prime})^2)}{(1-4(2kk^{\prime})^2)^{3/2}}+\frac{(k^\prime)^2-k^2}{(1-4(2kk^{\prime})^2)^{1/2}}-\frac{G_1}{(1-4(2kk^{\prime})^2)^{1/2}}\right).\end{displaymath}

\section{Some Preliminary Calculations}
\label{sec:a}

We will do some preliminary calculations before we move on to the main calculation in this document. Let $n=163$.

\begin{lemma}\label{lemma1} Denote $c(q):=c(k(q))$. Then $1728/c(e^{\pi i(2\tau+1)})=-j(\tau)$, where $j(\tau)$ is the Klein $j$-invariant.\end{lemma}

\begin{proof} From Chapter 2.1 of \cite{BorBro} one gets
\begin{equation}
\begin{split}
(2k(-q^2)k^\prime(-q^2))^2&=-4\left(\frac{\theta_2^2(q^2)\theta_3^2(q^2)}{\theta_4^4(q^2)}\right)^2\\
&=-\frac{1}{4}\left(\frac{(\theta_3^2(q)+\theta_4^2(q))(\theta_3^2(q)-\theta_4^2(q))}{\theta_3^2(q)\theta_4^2(q)}\right)^2\\
&=-\frac{1}{4}\left(\frac{\theta_2^4(q)}{\theta_3^2(q)\theta_4^2(q)}\right)^2.
\end{split}
\end{equation}
Then
\begin{equation}
\begin{split}
\frac{1728}{c(k(q))}&=-256\frac{\left(1+\theta_2^8(q)/(\theta_3^4(q)\theta_4^4(q))\right)^3}{\theta_2^8(q)/(\theta_3^4(q)\theta_4^4(q))}\\
&=-256\frac{\left(\theta_3^4(q)\theta_4^4(q)+\theta_2^8(q)\right)^3}{\theta_2^8(q)\theta_3^8(q)\theta_4^8(q)}\\
&=-32\frac{\left(\theta_3^8(q)+\theta_4^8(q)+\theta_2^8(q)\right)^3}{\theta_2^8(q)\theta_3^8(q)\theta_4^8(q)}
\end{split}
\end{equation}
Let $q=e^{\pi i\tau}$, then the function $1728/c(k(\exp(\pi i\tau)))$ is invariant under the transformation $\tau\mapsto\tau+1$ and $\tau\mapsto-1/\tau$. The $q$-expansion of $1728/c(k(q))=-1/q-744-\cdots$, and the result follows.
\end{proof}

A well-known result attributed to Hermite is that \begin{math}j\left(\frac{-1+\sqrt{-163}}{2}\right)=-640320^3\end{math}, which implies that \begin{equation}\label{cvalue}c(k(\exp(-\pi\sqrt{n})))=-\frac{12^3}{640320^3}=-\frac{1}{53360^3}.\end{equation}

We also remark that $2v(k)=\sqrt{1-c(k)}$, so
\begin{equation}\label{vvalue}
\begin{split}
2\sqrt{n}v(k)&=\sqrt{n(1-c(k))}\\
&=\sqrt{163}\sqrt{1+\frac{1728}{640320^3}}\\
&=\frac{12\times545140134}{640320^{3/2}}.\\
\end{split}
\end{equation}

The calculation of $2v(k)$ also implies that
\begin{equation}\label{Gvalue}
\begin{split}
G_0&=\frac{\sqrt{n}}{6}\left(2v(k)+\frac{(k^\prime)^2-k^2}{(1-4(2kk^{\prime})^2)^{1/2}}-\frac{G_1}{(1-4(2kk^{\prime})^2)^{1/2}}\right)\\
&=\frac{2\sqrt{n}v(k)}{6}+\frac{\sqrt{n}}{6}\frac{(2kk^\prime)^{1/3}}{(1-4(2kk^{\prime})^2)^{1/2}}\frac{(k^\prime)^2-k^2-G_1}{(2kk^\prime)^{1/3}}\\
&=\frac{2\sqrt{n}v(k)}{6}+\frac{\sqrt{n}}{3}\sqrt[6]{-\frac{c(k)}{1728}}\frac{(k^\prime)^2-k^2-G_1}{(2kk^\prime)^{1/3}}
\end{split}
\end{equation}

Combining (\ref{Gvalue}),(\ref{vvalue}),(\ref{cvalue}) and (\ref{chud}), one will be led to the statement to be proved later:

\begin{proposition}If $n=163$, then \begin{equation}\label{G2value}G_2:=\frac{\sqrt{n}}{3}\frac{(k^\prime)^2-k^2-G_1}{(2kk^\prime)^{1/3}}=-1448\end{equation}\end{proposition}
$\mathbf{Remark:}$ One can also get similar propositions for other imaginary quadratic number fields $\mathbb{Q}(\sqrt{-n})$ with class number 1. We list similar results in Table \ref{tab:g2val}. We also note that the values of $(2kk^\prime)^{1/3}$ at $\exp(-\pi\sqrt{n})$ are necessary in our calculation. We can extract the values of the cubic root of (Weber's) modular function at singular moduli from the special values of $j$-invariant and Lemma \ref{lemma1}. The (cubic) minimal polynomials of these values are listed in the Table \ref{tab:g2val} as well. The values of $(2kk^\prime)^{1/3}$ at singular moduli are the (unique) real roots of these polynomials.

\begin{table}[htbp]
\centering
 \caption{\label{tab:g2val}Values of $G_2$ and $(2kk^\prime)^{1/3}$ at singular moduli}
 \begin{tabular}{lcc}
  \toprule
  $n$ & $G_2$ &$\mathrm{min.\,pol.\,of\,}(2kk^\prime(\exp(-\pi\sqrt{n})))^{1/3}$\\
  \midrule
  19 & -4 & $2 x^3 + 4 x^2 + 4 x - 1$\\
  43 & -24 & $2 x^3 - 8 x^2 + 16 x - 1$\\
  67 & -76 & $2 x^3 + 12 x^2 + 36 x - 1$\\
  163 & -1448 & $2 x^3 + 40 x^2 + 400 x - 1$\\
  \bottomrule
 \end{tabular}
\end{table}

We need another lemma before we close the preliminary calculations.

\begin{lemma}\label{lemma2}Let $q_0=\exp({-\pi/\sqrt{n}})$. Then \begin{equation}\label{thetatrans}\left(\sum_{u,v\in\mathbb{Z}}q_0^{2u^2+2uv+(n+1)v^2/2}\right)^2=\frac{\theta_3^4(q_0)(1+2kk^\prime(q_0))}{2\sqrt{n}}.\end{equation}\end{lemma}

\begin{proof}
We denote the left-hand of (\ref{thetatrans}) by $Q_0.$
\begin{equation}
\begin{split}
Q_0&=\left(\theta_3(q_0^2)+\theta_2(q_0^2)\theta_2(q_0^{2n})\right)^2\\
&=\frac{\theta_3^2(q_0)\theta_3^2(q_0^{n})}{4}\left(\sqrt{(1+k(q_0))(1+k(q_0^n))}+\sqrt{(1-k^\prime(q_0))(1-k^\prime(q_0^n))}\right)^2\\
&=\frac{\theta_3^2(q_0)\theta_3^2(q_0^{n})}{4}\left(2+2kk^\prime(q_0)+2\sqrt{(1-(k(q_0))^2)(1-(k^\prime(q_0))^2}\right)\\
&=\frac{\theta_3^4(q_0)}{2\sqrt{n}}(1+2kk^\prime(q_0)).
\end{split}
\end{equation}
\end{proof}

\section{Construction of modular forms on $\Gamma_0(p)$}
\label{sec:b}

The main task of our calculation is to evaluate the constant $G_2$ in (\ref{G2value}). The most difficult part in the evaluation still lies in the calculation of the constant $G_1$ defined by (\ref{ramfuncs}). The main ingredient in our calculation originates from the work of E. Hecke, M. Eichler and A. Pizer on the construction of modular forms from quaternion algebra $A(p)$ ramified at a prime $p$ and $\infty$. A detailed theory of quaternion algebra is not our main concern, and we refer our readers to the paper \cite{Pizer} of A. Pizer.

Consider an integral even lattice $\Lambda$ of rank $4$ in the Euclidean space $\mathbb{R}_4$. The theta series associated to the lattice $\Lambda$ is defined to be $\theta_{\Lambda}(\tau)=\sum_{x\in\Lambda}q^{\lVert x\rVert^2},q=e^{\pi i\tau}$. Let $l^2$ be the determinant of $\Lambda$, $\Lambda^{*}$ be its dual lattice. If $\Lambda$ is isometric to $l\Lambda^{*}$, then we call the integral even lattice a self-dual modular lattice. The name originates from the theorem below\cite[Prop. 3.3]{Pizer1}\cite[Theorem 1]{Quebbemann}:

\begin{theorem}\label{modorigin}Let the determinant of a self-dual modular lattice $\Lambda$ of rank 4 be the square of an integer $l$. Then the theta function $\theta_{\Lambda}(\tau)$ associated to the modular lattice $\Lambda$ is a modular form of weight $2$ on $\Gamma_0(l)$. The modular form $\theta_{\Lambda}(\tau)$ is also an eigenfunction of the Fricke involution $w:z\mapsto-1/(lz)$, with eigenvalue -1, i.e., $\theta_{\Lambda}(-1/(l\tau))=-l\tau^2\theta_{\Lambda}(\tau)$.  \end{theorem}

Since the series $nP(q^n)-P(q),q=e^{\pi i\tau}$ is a modular form of weight $2$ with Fricke eigenvalue $-1$ as well, it is reasonable to represent this $q$-series by the linear combination of theta series attached to modular lattices with discriminant $n^2$. One will inevitably wonder how one can construct enough many modular lattices with the given discriminant $n^2$, and it is where the theory of quaternion algebra enters.

We limit us to the case $n=p$, $p$ a prime with $p\equiv3\pmod 4$ for the sake of conciseness. The quaternion algebra $A(p)$ is defined to be the (unique) central simple algebra over $\mathbb{Q}$ of dimension $4$ ramified at $p$ and $\infty$. In other words, the quaternion algebra is a 4-dimensional linear space over $\mathbb{Q}$ with a basis $1,i,j,k$ on which one can define multiplication: $i^2=-1,j^2=-p,ij=-ji=k$. The conjugate of an element $w=a+bi+cj+dk,a,b,c,d\in\mathbb{Q}$ is defined to be $\bar{w}=a-bi-cj-dk$. The trace $\mathrm{tr}(w)$ and the norm $N(w)$ are defined to be $\mathrm{tr}(w)=w+\bar{w}$ and $N(w)=w\bar{w}$. One can also define an inner product on $A(p)$ with respect to the norm: $\langle x,y\rangle:=N(x+y)-N(x)-N(y)=\mathrm{tr}(x\bar{y})$.

Since $A(p)$ is the non-commutative analogue of quadratic field over $\mathbb{Q}$, one can define a maximal order (not necessarily unique) $\mathcal{O}$ and (left) $\mathcal{O}$-ideal classes on it. A $\mathbb{Z}$-basis of a certain maximal order $\mathfrak{O}$ (with unity) can be given as follows\cite[Prop. 5.1]{Pizer}:$(1+j)/2,(i+k)/2,j,k$. The norm $N(I)$ of an $\mathfrak{O}$-ideal $I$ is defined to be the greatest common divisor of the norms for all elements in the ideal $I$. We claim that an $\mathfrak{O}$-ideal $I$ is a self-dual modular lattice when we equip the lattice with the inner product $\langle x,y\rangle/N(I),x,y\in I$.

\begin{theorem}\label{mainfricke}The theta function attached to an $\mathfrak{O}$-ideal I in the quaternion algebra $A(p)$ is defined to be \begin{displaymath}\theta_I(\tau)=\sum_{x\in I}q^{{2N(x)}/{N(I)}},q=e^{\pi i\tau}.\end{displaymath}Then $\theta_I(\tau)$ is a modular form of weight 2 on $\Gamma_0(p)$ with Fricke eigenvalue $-1$.\end{theorem}

\begin{proof}A classical result of theta functions associated to quaternary quadratic forms \cite[p. 106]{Pizer1} asserts that $\theta_I(\tau)$ is a modular form of weight 2. Notice that the maximal order $\mathfrak{O}$ is a self-dual modular lattice with respect to $\langle x,y\rangle$, then the second part of the theorem follows from Lemma 3.7 of \cite{Pizer1}.\end{proof}

A theorem attributed to Eichler (see \cite[Theorem 1.12]{Pizer}) asserts that the number of distinct $\mathfrak{O}$-ideal classes is finite. There are more than one algorithm to determine the basis of ideals in all ideal classes(e.g., \cite{Pizer1}\cite{KValg}). The algorithm in \cite{KValg} is adopted in the computer algebra system MAGMA, which is used in our calculation to generate modular theta functions attached to $\mathfrak{O}$-ideals.

We know from \cite[Prop. 2.17]{Pizer} that two isomorphic $\mathfrak{O}$-ideals have the same theta functions. So the number of different theta functions are upper bounded by the type number $T(p)$ of the quaternion algebra $A(p)$, say, the number of ideal classes that are not isomorphic. An intriguing question is whether all these theta functions are linear independent(see Remark 2.16 of \cite{Pizer}). We would like to go one step further to investigate whether these theta functions form a $\mathbf{basis}$ of weight 2 modular forms on $\Gamma_0(p)$ with Fricke eigenvalue $-1$. The investigation is based on a fact (probably known to Max Deuring) that the type number $T(p)$ coincides with the dimension of weight 2 modular forms on $\Gamma_0(p)$ with Fricke eigenvalue $-1$. They are given by (see \cite{Cohn}\cite[Theorem 6]{Quebbemann} for the explicit dimension formula(a result known to Robert Fricke), and see \cite[p. 93, Theorem B]{Pizer2} for the explicit type number formula):

\begin{equation}\label{typenum}
\begin{split}
T(p)=\frac{1+g}{2}+2^{-t_p-1}h(\sqrt{-p}),\\
t_p=\begin{cases}
1,&p\equiv1\pmod4,\\
-1,&p\equiv3\pmod8,\\
0,&p\equiv7\pmod8.\\
\end{cases}
\end{split}
\end{equation}
$g$ is the genus of $X_0(p)$:

\begin{equation}
g=\begin{cases}
\lfloor\frac{p+1}{12}\rfloor,&p\not\equiv1\pmod {12},\\
\lfloor\frac{p+1}{12}\rfloor-1,&p\equiv1\pmod {12},\\
\end{cases}
\end{equation}
$h(\sqrt{-p})$ is the class number of the imaginary quadratic field $\mathbb{Q}(\sqrt{-p})$. Then the linear independence of $T(p)$ theta functions suggests that they form a basis of the subspace of modular forms aforementioned. Numerical computations with MAGMA verify the conjecture for $p<227$, while for $p=227$ the 15 theta functions span a $13$-dimensional subspace.

\begin{table}
\caption{\label{tab:apmtrx}(Reduced) quadratic forms for ideals in $A(p),p=163$}
\begin{tabular}{cc}
$M_{I_1}=\begin{pmatrix} 2 & 0 & 0 & 1 \\ 0 & 2 & 1 & 0\\ 0 & 1 & 82 & 0\\ 1 & 0 & 0 & 82 \end{pmatrix}$ & $M_{I_2} = \begin{pmatrix} 4 & 0 & 1 & 2\\ 0 & 4 & 2 & 1\\ 1 & 2 & 42 & 1\\ 2 & 1 & 1 & 42 \end{pmatrix}$\\\\
$M_{I_3}=\begin{pmatrix} 6 & 0 & 1 & -2\\ 0 & 6 & 2 & 1\\ 1 & 2 & 28 & 0\\ -2 & 1 & 0 & 28 \end{pmatrix}$ & $M_{I_4} = \begin{pmatrix} 8 & 0 & 2 & -3\\ 0 & 8 & 3 & 2\\ 2 & 3 & 22 & 0\\ -3 & 2 & 0 & 22 \end{pmatrix}$\\\\
$M_{I_5}=\begin{pmatrix} 10 & 0 & 1 & -4\\ 0 & 10 & 4 & 1\\ 1 & 4 & 18 & 0\\ -4 & 1 & 0 & 18 \end{pmatrix}$ & $M_{I_6} = \begin{pmatrix} 12 & 0 & 1 & -2\\ 0 & 12 & 2 & 1\\ 1 & 2 & 14 & 0\\ -2 & 1 & 0 & 14 \end{pmatrix}$\\\\
$M_{I_7}=\begin{pmatrix} 12 & 0 & 2 & -5\\ 0 & 12 & 5 & 2\\ 2 & 5 & 16 & 0\\ -5 & 2 & 0 & 16\end{pmatrix}$ & $M_{I_8} = \begin{pmatrix} 14 & 0 & 5 & -6\\ 0 & 14 & 6 & 5\\ 5 & 6 & 16 & 0\\ -6 & 5 & 0 & 16 \end{pmatrix}$\\\\
\end{tabular}
\end{table}

We return to our case $p=163$. The type number $T(p)$ can be calculated with (\ref{typenum}), which is equal to $8$. We remark that $\theta_I$ can be rewritten as a theta function associated to a (positive-definite) quaternary quadratic form $x^{T}Mx$, where $M$ is a $4\times 4$ positive-definite matrix with integral entries:

\begin{equation}\theta_I(\tau)=\sum_{x\in\mathbb{Z}^4}q^{x^{T}Mx},q=e^{\pi i\tau}.\end{equation}
One can use quaternion algebra package in MAGMA to calculate $8$ matrices $M_{I_i}$ for all $\mathfrak{O}$-ideal classes $I_{1},\cdots,I_{8}$, then verify the linear independence of corresponding $8$ theta series from the coefficients of $q$-expansions. We list those eight integral matrices in the Table \ref{tab:apmtrx}.

$\mathbf{Remark:}$ The motivation of Lemma \ref{lemma2} is now clear. The left-hand side of (\ref{thetatrans}) is the theta function associated to the maximal order $\mathfrak{O}$ in the quaternion algebra $A(p)$.

\section{Construction of modular equations from theta functions associated to quaternion ideals}
\label{sec:c}

Since theta functions $\theta_{I_1},\cdots,\theta_{I_8}$ form a basis of modular forms of weight 2 on $\Gamma_0(p)$ with Fricke eigenvalue $-1$, the modular form $nP(q^n)-P(q)$ is a linear combination of theta series $\theta_{I_1},\cdots,\theta_{I_8}$. One can easily verify that
\begin{equation}\label{thetaspan}nP(q^n)-P(q)=6\theta_{I_1}+12\theta_{I_2}+24\sum_{i=3}^{8}\theta_{I_i}\end{equation} with a few coefficients from $q$-expansions of theta series and Eisenstein series. One is tempted to evaluate $\theta_{I_i}(q)/\theta_{I_1}(q)$ at $q=\exp(-\pi/\sqrt{n})$ when one gives a quick glimpse at (\ref{thetaspan}),  (\ref{ramfuncs}), (\ref{Kexp}) and(\ref{thetatrans}), but it is by no means convenient to work with the ratios $\theta_{I_i}(q)/\theta_{I_1}(q)$ directly. We would like to elaborate the method that Mazur and Swinnerton-Dyer used to construct models of $X_{0}^{+}(p)$ in their paper \cite{MaPS} for our case $n=163$(they construct a model for $n=37$ only).

$\mathbf{Remark:}$ It is technically convenient to work on $X_{0}^{+}(p)$ than $X_{0}(p)$, since the genus of $X_{0}^{+}(p)$ is generally half of the genus of $X_{0}(p)$. The genus for $X_{0}^{+}(p)$ is $6$, while the genus for $X_{0}(p)$ is $13$.

Theta functions arising from quaternion ideals can be used to construct meromorphic functions on $X_{0}^{+}(p)$: the ratios of the linear combinations of theta functions are meromorphic functions on the Riemann surface $X_{0}^{+}(p)$, and one can construct explicit models of modular curve $X_{0}^{+}(p)$ with these meromorphic functions. In order to construct simpler models for the curve $X_{0}^{+}(p)$, one should rather choose meromorphic functions with smaller degrees, i.e., meromorphic functions with fewer poles. One has to choose appropriate linear combination of theta functions with lowest number of zeros in the fundamental region of $\Gamma_0^{+}(p)$($\infty$ excluded). We claim that the cusp form $\phi=(\theta_{I_6}-\theta_{I_7})/4=q^{14}-q^{16}-q^{18}+\cdots$ is the linear combination with the lowest number of zero in the fundamental region of $\Gamma_0^{+}(p)$. An holomorphic $1$-form on the Riemann surfaces $X_0(p)$ has $2g-2=2\times13-2=24$ zeros, so $\phi$ has $(24-2(7-1))/2=6$ zeros other than the cusps on the Riemann surfaces $X_0^{+}(p)$. Note that $(\theta_{I_7}-\theta_{I_8})/4=q^{12}+\cdots$, then
\begin{equation}f=\frac{\theta_{I_7}-\theta_{I_8}}{\theta_{I_6}-\theta_{I_7}}\end{equation}
is the meromorphic function with smallest degree on the Riemann surface $X_0^{+}(p)$. In order to evaluate $f(\tau)$ at $\tau=i/\sqrt{163}$, we need to determine an explicit model of $X_0^{+}(p)$. We choose another meromorphic function
\begin{equation}
\varphi(\tau)=p^2\frac{\eta^4(p\tau)}{\eta^4(\tau)}+\frac{\eta^4(\tau)}{\eta^4(p\tau)}
\end{equation}
which is a meromorphic function on $X_0^{+}(p)$ when $p\equiv 7\pmod{12}$. $\varphi(\tau)$ is everywhere holomorphic on $X_0^{+}(p)$ except a pole at $\infty$ with order $27$. On the compact Riemann surface, two meromorphic function $\phi$(degree $=7$) and $\varphi$(degree $=27$) must satisfy an algebraic relation:

\begin{equation}\label{modeq}
\sum_{0\leq i\leq7}\varphi^{i}y_i(f)=0,
\end{equation}
while $y_i(t)$ are polynomials of $t$ with degree $\leq 27$.  With the first $200$ coefficients of the $q$-expansion of $f$ and $\varphi$, one can find the explicit expressions of the polynomials $y_i$.

We also define $g_i,i=1,\cdots,5$ as

\begin{equation}g_i(\tau)=\frac{\theta_{I_i}-\theta_{I_{i+1}}}{\theta_{I_6}-\theta_{I_7}}\end{equation}
and
\begin{equation}g_6(\tau)=\frac{4\theta_{I_6}}{\theta_{I_6}-\theta_{I_7}}\end{equation}
respectively. Since they are meromorphic functions on $X_0^{+}(p)$, one can also get similar modular equation as (\ref{modeq}) for $f$ and $g_i$. Once we get the value of $f$ at $\tau=i/\sqrt{163}$, we can use the modular equations for $f$ and $g_i$ to evaluate $g_i$ at $\tau=i/\sqrt{163}$.

\begin{theorem}Let $u$ be the unique real root of $2 x^3 + 40 x^2 + 400 x - 1=0$. We claim that $f(i/\sqrt{163})=(64u^2+1372u+11680)/4389$.\end{theorem}
\begin{proof}Let $f(i/\sqrt{163})=u^\prime$. We note that $\phi(i/\sqrt{163})=2\times 163$. One can plug this value in (\ref{mainmodeqst})-(\ref{mainmodeqed}) and one can get that $u^\prime$ is a root of a 27th degree polynomial
\begin{equation}
\begin{split}
(-&160 - 512 x - 400 x^2 + 231 x^3)\cdot(-44044178 - 500041768 x\\
- &1716358972 x^2 - 1771347457 x^3 + 949502158 x^4 + 1855221822 x^5\\
-  &251575929 x^6 - 706717664 x^7 + 79493657 x^8 + 121939618 x^9 \\
-  &18956160 x^{10} - 7891968 x^{11} + 1622016 x^{12})^2
\end{split}
\end{equation}
Numerical computation determines that $u^{\prime}$ is the root of the polynomial $-160 - 512 x - 400 x^2 + 231 x^3 = 0$, and the result follows.
\end{proof}

$\mathbf{Remark:}$ The zeros of $\phi(\tau)$ are very likely to be imaginary quadratic irrationals. $\phi(\tau)$ vanishes where $f(\tau)$ has a pole. From (\ref{mainmodeqst})-(\ref{mainmodeqed}) one can get that the value of $\varphi(\tau)$ at a pole of $f(\tau)$ is a root of the polynomial $-x^6 + 177 x^5 - 2442 x^4 - 1069320 x^3 + 49392000 x^2 + 989898000 x - 65739380000 = 0$. All the six roots of this polynomial are integers: $x_1=-70,x_2=-37,x_3=74,x_4=x_5=x_6=70$. One can conjecture that the zeros of $\phi(\tau)$ are elliptic points of $\Gamma_0(163)$ and some $Heegner\,points$ $\omega$ on $\Gamma_0(163)$, i.e., $\omega$ and $163\omega$ are imaginary quadratic irrationals whose minimal polynomials have the same discriminant $D$. Numerical computations suggests that $D=-7$(for $x_3=74$), $D=-11$(for $x_1=-70$), $D=-27$(for $x_2=-37$), $D=-44$(for $x_4=x_5=x_6=70$).

With the same method one can obtain all the values of modular function $g_1,\cdots,g_6$ at $\tau=i/\sqrt{163}$. We will list them in the Table \ref{tab:modval1} below.

\begin{table}[htbp]
\centering
 \caption{\label{tab:modval1}Values of $g_i$ and $f$ at singular moduli $i/\sqrt{163}$}
 \begin{tabular}{lc}
  \toprule
  $n$ & $\mathrm{min.\,pol.\,of\,values}$\\
  \midrule
  $f$ & $231 x^3- 400 x^2- 512 x -160 $\\
  $g_1$ & $231 x^3 - 36436 x^2 - 148896 x -211600$\\
  $g_2$ & $231 x^3  - 11843 x^2 - 19899 x -11153 $\\
  $g_3$ & $3 x^3  - 80 x^2  + 216 x -176$\\
  $g_4$ & $11 x^3 - 117 x^2 - 7 x - 335$\\
  $g_5$ & $11 x^3 - 19 x^2 - 63 x -81 $\\
  $g_6$ & $77 x^3 - 22020 x^2 + 86760 x -116964$\\
  \bottomrule
 \end{tabular}
\end{table}

One can extract the values of $\theta_{I_i}/\theta_{I_1}$ at $\tau=i/\sqrt{163}$ from Table \ref{tab:modval1}. Let $u$ still be the unique real root of $2 x^3 + 40 x^2 + 400 x - 1=0$.

\begin{table}[htbp]
\centering
 \caption{\label{tab:thetaratio}Values of $\theta_{I_i}/\theta_{I_1}$ at singular moduli $i/\sqrt{163}$}
 \begin{tabular}{lc}
  \toprule
  $n$ & $\mathrm{values\,of\,modular\,functions}$\\
  \midrule
  $\theta_{I_2}/\theta_{I_1}$ & $ (51856u^2 + 1015272u + 10681241) /21360009$\\
  $\theta_{I_3}/\theta_{I_1}$ & $ (644u^2+12408u+59401) /176529$\\
  $\theta_{I_4}/\theta_{I_1}$ & $(-124u^2+ 2592u +26905)/102201$\\
  $\theta_{I_5}/\theta_{I_1}$ & $(448u^2 -68928u+1634123) /7120003$\\
  $\theta_{I_6}/\theta_{I_1}$ & $ (868u^2+17064u+ 81593)/ 374737$\\
  $\theta_{I_7}/\theta_{I_1}$ & $ (-3792u^2+77288u+1528883)/7120003 $\\
  $\theta_{I_8}/\theta_{I_1}$ & $ (400u^2+19864u+133633)/ 647273$\\
  \bottomrule
 \end{tabular}
\end{table}

Combining Table \ref{tab:thetaratio}, (\ref{thetaspan}), (\ref{thetatrans}), we have

\begin{equation}
\begin{split}
G_2=&\frac{\sqrt{163}}{3u}\left(\sqrt{1-u^6}-\frac{1044640 u^2+ 33631776 u + 336327974}{7120003}\frac{1+u^3}{2\sqrt{163}}\right)\\
=&\frac{\sqrt{163}}{3u}\left(\sqrt{1-u^6}-\frac{40 u^2+ 748992 u + 8003}{209}\frac{1}{2\sqrt{163}}\right)\\
=&\frac{1}{3u}\left(\frac{40 u^2-1066800 u + 8003}{418}-\frac{40 u^2+ 748992 u + 8003}{418}\right)\\
=&-1448,
\end{split}
\end{equation}
and we are done.

$\mathbf{Remark:}$ We also note that the explicit models of $X_0^{+}(p)$ generated by these theta functions have comparatively small coefficients(with respect to the modular equation from $j(\tau),j(p\tau)$).

\section*{Acknowledgments}

The author would like to thank Prof. Heng Huat Chan and Dr. Jes\'{u}s Guillera for their valuable comments and suggestions on this manuscript.

\section{Appendix: Explicit modular equations for $n=163$}
\label{sec:d}
1. Modular equation for $f,\varphi$.
\begin{equation}\sum_{0\leq i\leq7}\varphi^{i}y_i(f)=0, y_7(x)=1;\end{equation}
\begin{equation}\label{mainmodeqst}
\begin{split}
y_0(x)=-&65739380000 x^{27} - 2738210160000 x^{26} - 12036747296000 x^{25} + 350898012641680 x^{24} +\\
&644683749050032 x^{23} - 10966748735578104 x^{22} - 18158847588921566 x^{21} +\\
 &170490800547851734 x^{20} + 279490703499480027 x^{19} - 1547392770518054918 x^{18}\\
 - &2324369674809215120 x^{17} + 8373387065260480638 x^{16} + 9804924454337019041 x^{15}\\
 - &26502508074913920502 x^{14} - 10497020053589572336 x^{13} + 53835578718190847240 x^{12}\\
 - &75578213019853201224 x^{11} - 133156026485879678504 x^{10} + 280929284170405966525 x^9 +\\
 &399760342254330512112 x^8 - 174258821253085699386 x^7 - 484048676226560027262 x^6 \\
 - &310839845578967657160 x^5 - 175173341489755599608 x^4 - 131300350040887496670 x^3 \\
 - &79481280928275478532 x^2 -28730725453703354936 x - 4772634032798114224
\end{split}
\end{equation}
\begin{equation}
\begin{split}
y_1(x)=&989898000 x^{27} + 43205926400 x^{26} + 399844120480 x^{25} - 12587240164160 x^{24}\\
- &10184364915352 x^{23} + 395645402631900 x^{22} + 367059214422379 x^{21}\\
- &5561964901073447 x^{20} -9118181309829979 x^{19} + 42710370269237119 x^{18} +\\
 &114274640716859171 x^{17} - 166503018578521639 x^{16} - 776610574395254280 x^{15} +\\
 &81149557481589740 x^{14} + 2881407475746413354 x^{13} + 1962600391132979131 x^{12}\\
- &5020300210556940676 x^{11} - 6994509469967819078 x^{10} + 1097875709711319396 x^9 +\\
&6506577948689732176 x^8 + 4267061979747189021 x^7 + 2911911004380789535 x^6 +\\
 &3787737947744127764 x^5 + 2881661414439000060 x^4 + 953258426607224943 x^3 +\\
 &29064250388446420 x^2 - 59701340867325260 x - 12096640200441408
\end{split}
\end{equation}
\begin{equation}
\begin{split}
y_2(x)=&49392000 x^{27} + 722397760 x^{26} - 3120028744 x^{25} + 143913481408 x^{24}\\
- &215801003088 x^{23} - 5277017188846 x^{22} + 5475896831338 x^{21} +\\
 &77033079426150 x^{20} - 21787659599462 x^{19} - 603852242182902 x^{18}\\
 - &386154558749083 x^{17} + 2490625204583624 x^{16} + 4744601996931698 x^{15}\\
 - &2674261902998852 x^{14} - 20711180387803771 x^{13} - 19844065558430734 x^{12} +\\
   &31007126937762651 x^{11} + 77327422584349144 x^{10} + 32251135025970508 x^9\\
 - &60305321091074134 x^8 - 95426598503401445 x^7 - 70323025560586048 x^6\\
 - &37098247964094146 x^5 - 17532083378121022 x^4 - 8015716918858076 x^3\\
 - &2691801414175607 x^2 - 403787531063906 x + 9967038795928\\
\end{split}
\end{equation}

\begin{equation}
\begin{split}
y_3(x)=-&1069320 x^{27} - 6967264 x^{26} - 72976364 x^{25} - 392365572 x^{24} +\\
 &4981341946 x^{23} + 23566594429 x^{22} - 108736773901 x^{21}\\
 -&450804990162 x^{20} + 1130037412277 x^{19} + 4573657934698 x^{18}\\
 -&5896700266436 x^{17} - 27327276404752 x^{16} + 11885497238541 x^{15} +\\
 &94379383455584 x^{14} + 17779182785984 x^{13} - 161825109682966 x^{12}\\
 -&98999957129288 x^{11} + 70187491324856 x^{10} + 11078698296477 x^9\\
 -&8656298967255 x^8 + 205274987829483 x^7 + 309761602483799 x^6 +\\
 &174843303301573 x^5 + 30924523734371 x^4 - 9793867731731 x^3\\
 -&4792460794049 x^2 - 380531544897 x + 87325254748
\end{split}
\end{equation}

\begin{equation}
\begin{split}
y_4(x)=-&2442 x^{27} - 157256 x^{26} + 2016206 x^{25} -12421 x^{24}\\
       -&55529427 x^{23}+ 54733067 x^{22} + 663922703 x^{21}\\
       -&535757547 x^{20} -4537814913 x^{19} + 246488234 x^{18} +\\
        &17991502301 x^{17} + 22176675934 x^{16} - 25840279853 x^{15}\\
       -&134996639944 x^{14} -107303225295 x^{13} + 291966282881 x^{12} +\\
 &562317988923 x^{11} + 4942976887 x^{10} - 784288390331 x^9\\
 -&740681694992 x^8 - 56281263902 x^7 + 354630493854 x^6 +\\
 &264213980939 x^5 + 93104310974 x^4 + 42157420946 x^3 +\\
 &24377740590 x^2 + 6719505200 x + 448428557
\end{split}
\end{equation}

\begin{equation}
\begin{split}
-y_5(x)=-&177 x^{27} - 1112 x^{26} + 16055 x^{25} + 18643 x^{24} - 408519 x^{23} +\\
 &69835 x^{22} + 4879873 x^{21} - 3134569  x^{20} - 33727084 x^{19} +\\
 &25992439  x^{18} + 147866336 x^{17} - 102941457 x^{16} - 419983431 x^{15} +\\
 &208546102 x^{14} + 726517006 x^{13} - 205544592  x^{12} - 593797855 x^{11} +\\
 &238405166 x^{10} - 12729530  x^9 - 686637595 x^8 + 11353904 x^7 +\\
 &740193815 x^6 + 371747255 x^5 - 158758934 x^4 - 215326942 x^3 \\
 -&82627632 x^2 - 13535351 x - 395279
\end{split}
\end{equation}

\begin{equation}\label{mainmodeqed}
\begin{split}
-y_6(x)=&x^{27} - 4 x^{26} - 25 x^{25} + 85 x^{24} + 373 x^{23} - 821 x^{22} - 3832 x^{21} + 4087 x^{20} + \\
 &26948 x^{19} - 5545 x^{18} - 128026 x^{17} - 56515 x^{16} + 400320 x^{15} +\\
 &395074 x^{14} - 772939 x^{13} - 1234453 x^{12} + 763610 x^{11} + 2177559 x^{10} - 18589 x^9\\
 -&2179861 x^8 - 748983 x^7 + 1126371 x^6 + 692065 x^5 - 201531 x^4\\
 -&220729 x^3 - 25052 x^2 + 11164 x + 1641\\
\end{split}
\end{equation}

2. Modular equation for $f,g_6$.
\begin{equation}\sum_{0\leq i\leq7}g_6^{i}y_i(f)=0, y_7(x)=-1;\end{equation}

\begin{equation}
\begin{split}
y_0(x)= -&48 x^{11} - 744 x^{10} - 3567 x^9 - 4154 x^8 + 21923 x^7 + 104098 x^6 +\\
 &216253 x^5 + 267896 x^4 + 209668 x^3 + 103792 x^2 + 32960 x + 6400
\end{split}
\end{equation}

\begin{equation}
\begin{split}
y_1(x)= -&16 x^{12} - 264 x^{11} - 1361 x^{10} - 1421 x^9 + 7653 x^8 + 19982 x^7 -\\
 &13059 x^6 - 127054 x^5 - 246958 x^4 - 250985 x^3 - 148484 x^2 -\\
 &54848 x - 12800
\end{split}
\end{equation}

\begin{equation}
\begin{split}
y_2(x)= &52 x^{11} + 638 x^{10} + 2259 x^9 + 279 x^8 - 13226 x^7 - 22942 x^6 +\\
 &10437 x^5 + 82466 x^4 + 119979 x^3 + 86746 x^2 + 37812 x + 10896
\end{split}
\end{equation}

\begin{equation}
\begin{split}
y_3(x)= -&68 x^{10} - 576 x^9 - 1083 x^8 + 1878 x^7 + 8756 x^6 + 7848 x^5 -\\
 &10264 x^4 - 28900 x^3 - 26542 x^2 - 13888 x - 5144
\end{split}
\end{equation}

\begin{equation}
\begin{split}
y_4(x)= 42 x^9 + 209 x^8 + 87 x^7 - 1091 x^6 - 2163 x^5 - 283 x^4 + 3568 x^3 + 4479 x^2 + 2877 x + 1460
\end{split}
\end{equation}

\begin{equation}
\begin{split}
y_5(x)= -11 x^8 - 29 x^7 + 36 x^6 + 195 x^5 + 166 x^4 - 194 x^3 - 394 x^2 -
 320 x - 250
\end{split}
\end{equation}

\begin{equation}
\begin{split}
y_6(x)= x^7 + x^6 - 5 x^5 - 10 x^4 + 2 x^3 + 14 x^2 + 15 x + 24
\end{split}
\end{equation}

3. Modular equation for $f,g_5$.
\begin{equation}\sum_{0\leq i\leq7}g_5^{i}y_i(f)=0, y_7(x)=1;\end{equation}

\begin{equation}
\begin{split}
y_6(x)= - x^2 + 4 x + 8
\end{split}
\end{equation}

\begin{equation}
\begin{split}
y_5(x)= -4 x^3 + 36 x + 36
\end{split}
\end{equation}

\begin{equation}
\begin{split}
y_4(x)= -6 x^4 - 20 x^3 + 44 x^2 + 157 x + 100
\end{split}
\end{equation}

\begin{equation}
\begin{split}
y_3(x)= -6 x^5 - 43 x^4 - 24 x^3 + 232 x^2 + 408 x + 189
\end{split}
\end{equation}

\begin{equation}
\begin{split}
y_2(x)= -4 x^6 - 43 x^5 - 105 x^4 + 87 x^3 + 567 x^2 + 646 x + 232
\end{split}
\end{equation}

\begin{equation}
\begin{split}
y_1(x)= -12 x^6 - 75 x^5 - 81 x^4 + 255 x^3 + 687 x^2 + 588 x + 174
\end{split}
\end{equation}

\begin{equation}
\begin{split}
y_0(x)= -8 x^6 - 22 x^5 + 43 x^4 + 241 x^3 + 365 x^2 + 241 x + 60
\end{split}
\end{equation}

4. Modular equation for $f,g_4$.
\begin{equation}\sum_{0\leq i\leq7}g_4^{i}y_i(f)=0, y_7(x)=-1;\end{equation}

\begin{equation}
\begin{split}
y_6(x)= x^3 - 2 x - 2
\end{split}
\end{equation}

\begin{equation}
\begin{split}
y_5(x)= x^2 - 3 x - 8
\end{split}
\end{equation}

\begin{equation}
\begin{split}
y_4(x)= x^4 + 4 x^3 + 9 x^2 + 7 x - 2
\end{split}
\end{equation}

\begin{equation}
\begin{split}
y_3(x)= -2 x^6 - 10 x^5 - 17 x^4 + 3 x^3 + 20 x^2 + 3 x - 7
\end{split}
\end{equation}

\begin{equation}
\begin{split}
y_2(x)= x^7 + 13 x^6 + 34 x^5 + 3 x^4 - 61 x^3 - 21 x^2 + 35 x + 20
\end{split}
\end{equation}

\begin{equation}
\begin{split}
y_1(x)= -3 x^7 - 23 x^6 - 29 x^5 + 50 x^4 + 70 x^3 - 29 x^2 - 40 x
\end{split}
\end{equation}

\begin{equation}
\begin{split}
y_0(x)= 2 x^7 + 10 x^6 + 2 x^5 - 25 x^4 + x^3 + 20 x^2
\end{split}
\end{equation}

5. Modular equation for $f,g_3$.
\begin{equation}\sum_{0\leq i\leq7}g_3^{i}y_i(f)=0, y_7(x)=-1;\end{equation}

\begin{equation}
\begin{split}
y_6(x)= x^5 - 3 x^3 - 3 x^2 - 9 x - 5
\end{split}
\end{equation}

\begin{equation}
\begin{split}
y_5(x)= 10 x^6 + 4 x^5 - 35 x^4 - 41 x^3 - 45 x^2 - 36 x - 11
\end{split}
\end{equation}

\begin{equation}
\begin{split}
y_4(x)= 39 x^7 + 33 x^6 - 152 x^5 - 241 x^4 - 177 x^3 - 117 x^2 - 61 x -
 13
\end{split}
\end{equation}

\begin{equation}
\begin{split}
y_3(x)= 78 x^8 + 99 x^7 - 327 x^6 - 703 x^5 - 539 x^4 - 253 x^3 - 128 x^2 -
 51 x - 8
\end{split}
\end{equation}

\begin{equation}
\begin{split}
y_2(x)= 85 x^9 + 138 x^8 - 381 x^7 - 1063 x^6 - 999 x^5 - 457 x^4 -
 131 x^3 - 56 x^2 - 18 x - 2
\end{split}
\end{equation}

\begin{equation}
\begin{split}
y_1(x)= 48 x^{10} + 91 x^9 - 233 x^8 - 802 x^7 - 939 x^6 - 536 x^5 -
 111 x^4 + 8 x^3 - 5 x^2 - x
\end{split}
\end{equation}

\begin{equation}
\begin{split}
y_0(x)= 11 x^{11} + 23 x^{10} - 59 x^9 - 239 x^8 - 341 x^7 - 253 x^6 - 70 x^5 +
 21 x^4 + 9 x^3 - 2 x^2
\end{split}
\end{equation}

6. Modular equation for $f,g_2$.
\begin{equation}\sum_{0\leq i\leq7}g_2^{i}y_i(f)=0, y_7(x)=-1;\end{equation}

\begin{equation}
\begin{split}
y_6(x)= x^4 - 3 x^2 + 5 x + 6
\end{split}
\end{equation}

\begin{equation}
\begin{split}
y_5(x)= -7 x^5 - 4 x^4 + 20 x^3 + 4 x^2 - 19 x - 14
\end{split}
\end{equation}

\begin{equation}
\begin{split}
y_4(x)= 18 x^6 + 17 x^5 - 45 x^4 - 39 x^3 + 4 x^2 + 24 x + 18
\end{split}
\end{equation}

\begin{equation}
\begin{split}
y_3(x)= -20 x^7 - 16 x^6 + 42 x^5 + 23 x^4 + 6 x^3 + 4 x^2 - 18 x - 17
\end{split}
\end{equation}

\begin{equation}
\begin{split}
y_2(x)= 8 x^8 - 12 x^7 - 22 x^6 + 67 x^5 + 25 x^4 - 59 x^3 - 21 x^2 +
 18 x + 12
\end{split}
\end{equation}

\begin{equation}
\begin{split}
y_1(x)= 16 x^8 + 12 x^7 - 59 x^6 - 25 x^5 + 63 x^4 + 35 x^3 - 8 x^2 -
 14 x - 4
\end{split}
\end{equation}

\begin{equation}
\begin{split}
y_0(x)= 8 x^7 + 10 x^6 - 16 x^5 - 11 x^4 + 6 x^3 + 3 x^2 + 4 x
\end{split}
\end{equation}

7. Modular equation for $f,g_1$.
\begin{equation}\sum_{0\leq i\leq7}g_1^{i}y_i(f)=0, y_7(x)=1;\end{equation}

\begin{equation}
\begin{split}
y_6(x)= - x^6 - x^5 + 4 x^4 + 8 x^3 - 2 x^2 - 12 x - 1
\end{split}
\end{equation}

\begin{equation}
\begin{split}
y_5(x)= 4 x^7 + 9 x^6 - 9 x^5 - 39 x^4 - 10 x^3 + 35 x^2 + 16 x + 8
\end{split}
\end{equation}

\begin{equation}
\begin{split}
y_4(x)= -7 x^8 - 30 x^7 - 11 x^6 + 89 x^5 + 100 x^4 - 23 x^3 - 77 x^2 -
 52 x - 4
\end{split}
\end{equation}

\begin{equation}
\begin{split}
y_3(x)= 6 x^9 + 45 x^8 + 70 x^7 - 46 x^6 - 171 x^5 - 101 x^4 + 15 x^3 +
 71 x^2 + 48 x + 16
\end{split}
\end{equation}

\begin{equation}
\begin{split}
y_2(x)= - x^{10} - 27 x^9 - 86 x^8 - 68 x^7 + 31 x^6 + 74 x^5 + 72 x^4 +
 28 x^3 + 4 x^2
\end{split}
\end{equation}

\begin{equation}
\begin{split}
y_1(x)= -2 x^{11} + 26 x^9 + 50 x^8 + 48 x^7 + 36 x^6 + 12 x^5 + 4 x^4
\end{split}
\end{equation}

\begin{equation}
\begin{split}
y_0(x)= x^{12} + 4 x^{11} + 6 x^{10} + 6 x^9 + 5 x^8 + 2 x^7 + x^6
\end{split}
\end{equation}

\end{document}